\documentclass[11pt]{amsart}
\usepackage{a4wide}
\usepackage{hyperref}
\usepackage{amsthm,amssymb}
\usepackage[english]{babel}
\newcommand{\RR}{\mathbb{R}}
\newcommand{\CC}{\mathbb{C}}
\newcommand{\NN}{\mathbb{N}}
\newcommand{\ZZ}{\mathbb{Z}}

\newcommand{\PP}{\mathbb{P}}

\renewcommand{\varLambda}{\mathcal{V}}

\newcommand{\cA}{\mathcal{A}}
\newcommand{\cB}{\mathcal{B}}
\newcommand{\cC}{\mathcal{C}}
\newcommand{\cD}{\mathcal{D}}

\newcommand{\cF}{\mathcal{F}}
\newcommand{\cG}{\mathcal{G}}

\newcommand{\cL}{\mathcal{L}}

\newcommand{\cals}{\mathcal{F}}
\newcommand{\cala}{\mathcal{A}}
\newcommand{\calpM}{\mathcal{P}_0^B (M)}

\newcommand{\tr}{\mbox{tr}}

\newcommand{\dd}{\mathrm{d}}

\newcommand{\Tr}{{\mathop{\mathrm{Tr}}}}

\DeclareMathOperator{\supp}{\mathrm{supp}}

\newcommand{\Gammaf}{\Gamma_{\mbox{full}}}

\newcommand{\Gammat}{\Gamma_{\mbox{trans}}}

\newcommand{\be}{\begin{equation}}
\newcommand{\ee}{\end{equation}}
\newcommand{\bea}{\begin{eqnarray*}}
\newcommand{\eea}{\end{eqnarray*}}

\newtheorem{thm}{Theorem}
\newtheorem{lem}[thm]{Lemma}
\newtheorem{prp}[thm]{Proposition}
\newtheorem{cor}[thm]{Corollary}
\theoremstyle{definition}
\newtheorem{dfn}[thm]{Definition}
\newtheorem{ass}{Assumption}
\theoremstyle{remark}
\newtheorem{rem}[thm]{Remark}

\newcommand{\Hm}[1]{\leavevmode{\marginpar{\tiny%
$\hbox to 0mm{\hspace*{-0.5mm}$\leftarrow$\hss}%
\vcenter{\vrule depth 0.1mm height 0.1mm width \the\marginparwidth}%
\hbox to
0mm{\hss$\rightarrow$\hspace*{-0.5mm}}$\\\relax\raggedright #1}}}

\begin{document}
\title[Uniform existence of the IDS]{Uniform existence of the integrated density of states for random Schr\"odinger operators on metric graphs over $\ZZ^d$}

\author[M.~J.~Gruber]{Michael J.\ Gruber}
\author[D.~H.~Lenz]{Daniel H.\ Lenz}
\author[I.~Veseli\'c]{Ivan Veseli\'c}

\address[]{Fakult\"at f\"ur Mathematik, D-09107 Chemnitz, Germany}
\email{dlenz@mathematik.tu-chemnitz.de }
\urladdr{http://www.tu-chemnitz.de/mathematik/analysis/dlenz}
\urladdr{http://www.tu-chemnitz.de/~mjg/}
\urladdr{http://www.tu-chemnitz.de/mathematik/schroedinger/members.php}

\date{\today,  \jobname.tex}

 \subjclass[2000]{}

\begin{abstract}
We consider ergodic random Schr\"odinger operators on the metric graph
$\ZZ^d$ with random potentials and random boundary conditions taking
values in a finite set.  We show that normalized finite volume
eigenvalue counting functions converge to a limit uniformly in the
energy variable. This limit, the integrated density of states, can be
expressed by a closed Shubin-Pastur type trace formula.  It supports
the spectrum and its points of discontinuity are characterized by
existence of compactly supported eigenfunctions. Among other examples we
discuss random magnetic fields and percolation models.
\end{abstract}

\maketitle
\let\languagename\relax

\bigskip
\section{Introduction}

This paper is devoted to the spectral analysis of certain random ergodic operators
defined on metric graphs with a $\ZZ^d$-structure.
More precisely, we consider a Schr\"odinger operator with an ergodic random potential which takes values
in a finite set of potentials, with ergodic random boundary conditions from a finite set, 
and with ergodic random magnetic fields from a finite set. Such a model may be called random quantum graph.
A restriction of such an operator to a finite cube with selfadjoint boundary conditions
has purely discrete spectrum, hence it has a well defined eigenvalue counting function.
As the size of the cube tends to infinity the normalized counting functions converge
to the so-called \textit{integrated density of states} (IDS) or \textit{spectral distribution function}.
The corresponding result is well established for ergodic Schr\"odinger operators in the continuum (i.e.~on $L^2(\RR^2)$)
and on the lattice (i.e.~on $\ell^2(\ZZ^d)$), see for instance \cite{Pastur-80,KirschM-82c,Kirsch-89a,CarmonaL-90,PasturF-92}.
A modification of the proof to ergodic, random Schr\"odinger on metric graphs was presented in \cite[\S~6]{HelmV} 
and yields the existence of the IDS.
\smallskip

However, the abovementioned results concern \emph{pointwise} convergence of the distribution functions, or rather
\emph{weak convergence} of the corresponding measures. In the present paper we show that for the
class of random operators at hand the convergence holds uniformly in the energy variable.

Two features of the operators considered are helpful to derive such a strong statement about the convergence:
Firstly, since the potential and other data take only finitely many values frequencies of finite configuration patterns
exist in a nice way. Secondly, since the underlying geometric structure exhibits in many ways one-dimensional properties
the relevant perturbation operators turn out to have a uniformly bounded spectral shift function.
The latter feature can be understood as a finite-rank property in a generalized sense,
see Section \ref{s-FRP} for a more precise statement.

These two properties of the considered operators enable us to apply an abstract ergodic theorem for
Banach space valued random variables which was proven recently in \cite{LenzMV}. There the ergodic theorem
was already used to establish uniform convergence for a certain class of random Schr\"odinger operators
on $\ell^2(\ZZ^d)$. The situation for Schr\"odinger operators on metric graphs is somewhat more complicated
since the operators are unbounded and the IDS is not a probability distribution any more. Still,
quantum graph operators are technically much easier to handle than Schr\"odinger operators in truly
multi-dimensional space, i.~e.~on $L^d(\RR^2), d \ge 2$.

\bigskip

The paper is organized as follows. In Section \ref{s-MR} we define
precisely the operators we want to study and state the main results.
In Section \ref{s-FRP} we discuss finite rank perturbations of the
kind as they appear in our proofs.  In Section \ref{s-ET} we show that
the ergodic theorem of \cite{LenzMV} can be applied to the sequence of
spectral shift functions obtained from an exhaustion of the metric
graph by cubes.  This is then used in Section \ref{Proofs} to prove
our main results. In Section~\ref{magnetic} we show how magnetic fields on metric graphs
manifest themselves in the boundary conditions. Finally, in Section \ref{Site} we discuss
site- and edge-percolation  as an example for our results.

\section{Model and results}\label{s-MR}

We define a metric graph $G_d$ over $\ZZ^d$ in the following way. Let
$e_j$, $j=1,\ldots,d$ be the standard basis of the real $d$-dimensional space $\RR^d$.  The vertex
set $V_d$ of the graph consists of the points $\ZZ^d \subset \RR^d$ which have integer coefficients. Now, to
each vertex $x$ and $j\in\{1,\ldots,d\}$ we associate the edge $e=[x,x + e_j]$ with starting
point $s(e) = x$ and endpoint $r(e) =x+e_j$ given by
$$[x,x + e_j] :=\{x + t e_j : t \in]0,1[ \}.$$ Thus, each edge can
be canonically identified with the interval $]0,1[$.  This procedure may
seem to induce an orientation on our graph. However, it turns out that
all relevant quantities are independent of the choice of orientation.

The set of all edges is denoted by $E_d$. The union $V_d \cup\bigcup_{e\in E_d}e$ is
a closed subset of $\RR^d$, hence we can consider it as a topological
space with the metric inherited from $\RR^d$.  We denote this metric
subspace by $G_d$.  We will need to consider finite subgraphs $G$ of
$G_d$ as well. By a subgraph we mean a subset of the edges of $G_d$
together with all adjacent vertices. All functions we consider live on
the topological space $G_d$ or subgraphs of it. This shows, in
particular, that the choice of orientation at each edge is not
relevant.

\medskip

The operators we are interested in will be defined on the Hilbert space
\[
L^2(E_d):= \bigoplus_{e\in E_d } L^2(e)
\]
and their domains of definition will be subspaces of
\[
W^{2,2}(E_d):= \bigoplus_{e\in E_d } W^{2,2}(e),
\]
where $ W^{2,2}(e)$ is the usual Sobolev space of $L^2(e)$ functions
whose (weak) derivatives up to order two are in $L^2(e)$ as well. 
The restriction of $f\in W^{2,2}(E_d) $ to an edge $e$ is denoted by $f_e$.
  For an edge $e=[s(e),r(e)]=[x,x+e_j]$ and $g$ in $W^{2,2}(e)$ the
boundary values of $g$
$$ g(s(e)):=\lim_{t\searrow 0} g (x + t e_j), \;\:
g(r(e)):=\lim_{t\nearrow 1} g(s(e) + t e_j)$$  and the boundary
values of $g'$
\[
g'(s(e)):=\lim_{\epsilon \searrow 0} \frac{g(x+\epsilon e_j)-g(x)}{\epsilon} 
\quad \text{ and } \quad 
g'(r(e)):= \lim_{\epsilon \searrow 0} \frac{g(r(e))-g(r(e) -\epsilon e_j)}{- \epsilon}
\] 
exist by standard Sobolev type theorems. Note that we have
introduced a sign in the derivative at the endpoint of an edge. This
makes our definition of derivative canonical, i.e.\ independent of the
choice of orientation of the edge.
For  $f\in W^{2,2}(E_d)$ and each vertex $x$ we gather the boundary values
of $f_e (x)$ over all edges $e$ adjacent to $x$ in a vector
$f(x)$. Similarly, we gather the boundary values of $f_e' (x)$ over all
edges $e$ adjacent to $x$ in a vector $f' (x)$.

Given the boundary values of functions, we can now dicuss the concept
of boundary condition. Here we use material from \cite{KostrykinS-99b,Harmer-00}
to which we
refer for further details and proofs.  A single-vertex boundary condition at $x\in V$ 
is a choice of subspace $L_x$ of $\CC^{4d}$ with
dimension $2d$ such that
 $$ 
 \eta((v,v'),(w,w')):=\langle v',w\rangle - \langle v, w'\rangle
 $$ 
vanishes for all $(v,v'), (w,w')\in L_x$.  An $f\in W^{2,2}(E_d)$ is said to
satisfy the single-vertex boundary condition $L_x$ at $x$ if $(f(x),f'(x))$ belongs
to $L_x$.
A field  of single-vertex boundary conditions $ L:=\{L_x : x\in V_d\}$ will be called 
boundary condition. Given such a field, we obtain a selfadjoint realization $\Delta_L$ 
of the Laplacian $\Delta$ on $L^2 (E_d)$ by choosing the domain
 $$ 
 \cD(\Delta_L) :=\{f\in W^{2,2} (E_d) :\forall x : (f(x),f'(x))\in L_x\}.
 $$

Particularly relevant boundary conditions are Dirichlet boundary
conditions with subspace $L^D$ consisting of all those $(v,v')$ with $v=0$, Neumann conditions with
subspace $L^N$ consisting of all those $(v,v')$
with $v'=0$, and Kirchhoff (also known as free) boundary conditions $L^K$ consisting
of all $(v,v')$ with $v$ having all components equal and $v'$ having the sum over its components equal to $0$.

\begin{rem} 
More general types of boundary conditions are
conceivable. We will restrict ourselves to those just introduced. They may be
called graph local boundary conditions, as each boundary condition
involves values of $f$ at one vertex only.
\end{rem}

Everything discussed so far including existence of limits of functions at
the vertices and the notions of boundary condition extends in the
obvious way to subgraphs.  Moreover, for a subgraph $G$ of $G_d$ with edge set $E$, we
write $W^{2,2} (E):=\oplus_{e\in E} W^{2,2} (e)$. The number of edges of
a finite subgraph $G$ of $G_d$ is denoted by $|E|$.

\medskip

 In order to define random operators we need some further data
including a probability space $(\Omega,\PP)$ and an action of (a
subgroup of) the automorphism group of $G_d$ on $\Omega$ and maps $L$,
$V$ from $\Omega$ into the space of boundary conditions and potentials,
respectively.

For us two groups will be relevant, the full
automorphism group $\Gammaf$ and the group $\Gammat$ of translations
by $\ZZ^d$.  Note that $\Gammaf$ is generated by translations by
vectors in $\ZZ^d$ and a finite set of rotations.  We choose $\Gamma$
to be $\Gammaf$ or $\Gammat$ and assume that it acts ergodically on
$(\Omega,\PP)$ via measure preserving transformations. To simplify the
notation 
we identify $\gamma\in \Gamma$ with
the associated measure preserving transformation.

\bigskip

Let us describe the type of random operators we will consider in this paper:

 \begin{ass}\label{(S)}
Let $(\Omega,\PP)$ be a probability space and $\Gamma\in \{\Gammaf,\Gammat\}$ a group 
acting ergodically on $(\Omega,\PP)$. 
Let $\cB$ be a finite subset of
$L^\infty(0,1)$ and $\cL$ a finite set of boundary conditions. A random potential is a map
$$V : \Omega\longrightarrow \bigotimes_{e\in E_d} \cB \;\:\mbox{with}\;\:\;
V(\gamma (\omega))_{\gamma (e)} = V(\omega)_{e}$$
for all $\gamma\in \Gamma$ and $e\in
E_d$. A random boundary condition is a map 
$$ L : \Omega\longrightarrow \bigotimes_{v\in V_d} \cL,\;\:\mbox{with}\;\:\; L(\gamma (\omega))(\gamma (x)) = L(\omega)(x)$$
for all $\gamma\in \Gamma$ and $v\in V_d$.

A family of random operators $(H_\omega)$
on $L^2 (E_d)$ can be defined with domain of definition
$$ \cD(H_\omega ):=\{f\in W^{2,2} (E_d) : (f(x),f' (x)) \in L(\omega)(x)\;\:\mbox{ for all $x\in V_d$}\}$$
acting by
$$ (H_\omega  f)(e) := - f_e'' + V(\omega)_e f_e$$
for each edge $e$.    These are selfadjoint lower bounded operators.
 \end{ass}

We assume throughout the paper that Assumption \ref{(S)} holds, 
and for this reason do not repeat it in every lemma. Note however,
that in many statements we will only need part of the structure described in 
Assumption \ref{(S)}.

\begin{rem}  While $\Gammaf$ is not commutative it is a natural
object to deal with. In particular, let us note that the Laplacian
without potential and boundary conditions in all vertices identical to
Kirchhoff conditions is invariant under $\Gammaf$.
\end{rem}

We will need to consider restrictions of our operators to finite subgraphs.  These are finite subgraphs associated to finite subsets of $\ZZ$. The cardinality of a finite subset $Q$ of $\ZZ^d$ is denoted by $|Q|$.
Let  $$S:=\{\sum_{j=1}^d t_j e_j : 0< t_j < 1 \}.$$
  For a finite subset
$Q$ of $\ZZ^d$, we define the associated subgraph $G_Q$ of $G_d$ by
$$ G_Q := (Q + \bar S) \cap G_d.$$ 
The vertices and edges of $G_d$ contained in $G_Q$
are denoted by $V_Q$ and $E_Q$, respectively. Note that $V_Q \supset Q$. 
The set $V^i_Q$ of inner
vertices of $G_Q$ is then given by those vertices of $G_Q$ all of
whose adjacent edges (in $G_d$) are contained in $G_Q$.  The set of inner edges
$E^i_Q$ of $G_Q$ is given by those edges whose both endpoints are inner.
The vertices of $G_Q$ which are not inner are called boundary
vertices. The set of all boundary vertices is denoted by
$V^\partial_Q$.  Similarly, the set of edges which are not inner, is
denoted by $E^\partial_Q$.

\smallskip

The restriction $H_\omega^Q$  of the  random operator $H_\omega$  to  $G_Q$ has domain given by
\begin{align*}
\cD(H_\omega^Q):=
\{f \in W^{2,2} (E_Q)\mid & \forall x \in V^i_Q : (f(x),f'(x))\in L(\omega)(x), \\
 &\forall x \in V^\partial_Q : (f(x),f'(x)) \in L^D\}
\end{align*}
This operator is again selfadjoint, lower bounded, and has purely discrete spectrum. Let us enumerate the eigenvalues
of $H_\omega^Q$ in ascending order
$$\lambda_1(H_\omega^Q) < \lambda_2(H_\omega^Q) \le \lambda_3(H_\omega^Q) \le \dots$$
and counting multiplicities.
Then, the eigenvalue counting function $n_\omega^Q$ on $\RR$ defined by
\[
n_\omega^Q (\lambda) := \sharp \{n \in \NN \mid \lambda_n(H_\omega^Q) \le \lambda\}
\]
is monotone increasing and right continuous, i.e.~a distribution function,
which is associated to a pure point measure, $\mu_\omega^Q$.  Denote by
$$
N_\omega^Q(\lambda) := \frac{1}{ |E_Q|} n_\omega^Q (\lambda)
$$
the volume-scaled version of $n_\omega^Q(\lambda)$ and note that
$ |E_Q| = d |Q|$
as the edge to vertex ratio in the graph $(V_d,E_d)$ is equal to $d$.

 A sequence $(Q_l)_{l\in\mathbb{N}}$ of
finite subsets of $\ZZ^d$ is called a \emph{van Hove sequence} in
$\ZZ^{d}$ if  $
  \lim_{l  \to \infty} \frac{|V^\partial_{Q_l}|}{|Q_l|}=0. $
  For a finite subgraph $H$ of $G_d$ let
$\chi_{H }$ the multiplication operator by the characteristic function
of  $H$. Denote the trace on the operators on $L^2(E_d)$ by  $\Tr[\dots]$.
We can now state the main result of the paper.

\begin{thm}
\label{t-uniformIDS}
Let $Q$ be a finite subset of $\ZZ^d$.
Then, the function $N$ defined by
 \begin{equation} \label{e-IDS}
 N(\lambda):= \frac{1}{|E_Q|}\int_\Omega   \Tr\left[\chi_{G_Q} \chi_{]-\infty,\lambda]}(H_\omega) \right]  d \PP(\omega)
 \end{equation}
does not depend on the choice of $Q$,  is the distribution function of a measure $\mu$,  and
for any
van Hove sequence $(Q_l)$ in $\ZZ^d$
 \[
 \lim_{l \to\infty} \|N_\omega^{Q_l} -N\|_\infty =0
 \] 
for almost every $\omega \in \Omega$. In particular, for almost every
$\omega\in\Omega$, $N_\omega^{Q_l} (\lambda)$ converges as $l\to\infty$ pointwise to $N(\lambda)$ for every
$\lambda\in\RR$.
\end{thm}

\begin{rem}
(a) \
We chose to define the integrated density of states by a Shubin-Pastur formula 
and to prove that it coincides with the almost sure limit along van Hove sequences using an ergodic theorem.
Alternatively, one could take the point of view that the intuitively relevant 
objects are the normalised eigenvalue counting functions on finite graphs. Then the IDS
would be defined as their limit along a van Hove sequence and the trace formula \eqref{e-IDS} would be the result
of the theorem.

(b) \
For random Schr\"odinger operators with more general potentials and Kirchhoff boundary conditions on all vertices
the IDS was constructed in \cite{HelmV} as the pointwise almost everywhere limit of the sequence $N_\omega^{Q_l}$.

(c) \
In the case that $\Gamma$ is the full automorphism group $\Gammaf$ 
the equality
\begin{equation}\label{e-IDS-Gammaf}
 N(\lambda)= \frac{1}{|E(G)|}\int_\Omega   \Tr\left[\chi_{G} \chi_{]-\infty,\lambda]}(H_\omega) \right]  d \PP(\omega)
\end{equation}
holds for any finite subgraph $G$ with edge set $E(G)$.
To see this note that by linearity 
\begin{align*}
 \int_\Omega   \Tr\left[\chi_{S} \chi_{]-\infty,\lambda]}(H_\omega) \right]  d \PP(\omega)
&=
\sum_{j=1}^d \int_\Omega   \Tr\left[\chi_{[0,e_j]} \chi_{]-\infty,\lambda]}(H_\omega) \right]  d \PP(\omega)
\\&=
d \int_\Omega   \Tr\left[\chi_{[0,e_1]} \chi_{]-\infty,\lambda]}(H_\omega) \right]  d \PP(\omega).
\end{align*}
For the second equality we used that the $d$ terms in the sum are all equal, 
since the group $\Gammaf$ acts transitively on the edges.
For an arbitrary finite graph $G$ one can do a similar calculation, and formula 
\eqref{e-IDS-Gammaf} follows.
\end{rem}

While the definition of the IDS involves an ergodic theorem, there are other spectral features
of $H_\omega$ whose almost sure independence of $\omega$ uses only the ergodicity of the group action.
Prominent examples are the spectrum $\sigma(H_\omega)$ and its subsets 
$\sigma_{pp}(H_\omega)$, $\sigma_{sc}(H_\omega)$, $\sigma_{ac}(H_\omega)$, $\sigma_{disc}(H_\omega)$, $\sigma_{ess}(H_\omega)$ 
according to the spectral type.  In fact, by applying the general framework of \cite{LenzPV-02} we immediately infer the following theorem. 
\begin{thm} There exist subsets of the real line $\Sigma$, $\Sigma_{pp}$, $\Sigma_{sc}$, $\Sigma_{ac}$, $\Sigma_{disc}$, $\Sigma_{ess}$ 
and an $\Omega'\subset\Omega$ of full measure such that $\sigma(H_\omega)=\Sigma$ 
 and $\sigma_\bullet(H_\omega)=\Sigma_\bullet$ for all these spectral types $\bullet \in \{pp,sc,ac,disc,ess\}$
and all $\omega \in\Omega'$.
\end{thm}

Our goal here  is to establish a relation between the almost sure spectrum of $H_\omega$ and its	 IDS.
More precisely, in two corollaries to Theorem \ref{t-uniformIDS}
we relate the topological support $\supp \mu$ of  $\mu$ and the set $S_p (\mu):=\{\lambda\in \RR: \mu(\{\lambda\})>0\}$ of atoms of  $\mu$ to the spectrum of 
$H_\omega$. Note that the set of discontinuities of the IDS is precisely  $S_p (\mu)$.

\begin{cor}\label{spectrum}
$\Sigma$ equals the topological support $\supp \mu$ of $\mu$.
\end{cor}

As usual an $f\in W^{2,2} (E_d)$ is said to be compactly supported if $f_e\equiv 0$ for all but finitely many edges $e$.

\begin{cor}\label{jump}
Denote by $\Sigma_{cmp}$ the set of energies $\lambda\in \RR$
such that there exists almost surely a compactly supported $L^2$-eigenfunction $f_\omega$ with
$H_\omega f_\omega = \lambda f_\omega$. Then
\begin{equation} \label{e-jumps}
S_p(\mu) = \Sigma_{cmp}.
\end{equation}
\end{cor}

\begin{rem} 
(a) \
We will prove equation \eqref{e-jumps} using the following ingredients: Theorem \ref{t-uniformIDS} and 
the reasoning presented in \cite{KlassertLS-03}. Note, however, that both of the sets in 
\eqref{e-jumps} are given in terms of the infinite-volume operator $H_\omega$ 
rather than by its restrictions to finite subgraphs. Indeed, there exists a proof of \eqref{e-jumps}
which makes use neither of finite-graph operators, respectively their eigenvalue counting functions,
nor of an ergodic theorem: 
Once one has established the invariance property of the IDS, as done in e.g.~\cite{LenzPV-02}, equation \eqref{e-jumps} follows from the mere existence of a van Hove sequence.
This strategy is implemented in the proof of Theorem 2.3~(i) of \cite{Veselic-05b} for operators on combinatorial graphs
and can also be adapted to our situation.
Note however, that there are certain situations where it is rather natural to relate the (size of the) jumps of the IDS
to those of its finite-graph analogues, see for instance \cite{DodziukLMSY-03,Elek-06,LenzV}.

(b) \ 
There are many examples where the IDS has discontinuities. 
Indeed, the free Laplacian, i.e.~the Schr\"odinger operator with identically vanishing potential,
with Dirichlet, Neumann, or Kirchhoff boundary conditions exhibits compactly supported eigenfunctions. 
The relevance of Theorem \ref{t-uniformIDS} is that one has uniform
convergence all the same.  

(c) \
If the randomness entering the potential of the operator is sufficiently strong it is natural
to expect a smoothing effect on the IDS. In fact, in \cite{HelmV} for a class of alloy-type random potentials
the Lipschitz-continuity of the IDS was established. 
In \cite{GruberV} we  show for a different class of random potentials how one can 
estimate the modulus of continuity of the IDS.
\end{rem}

\section{Finite rank perturbations}\label{s-FRP}
We will have to deal with changes  in the spectral counting function due to
changes of boundary conditions and changes of
potentials. In particular, we will have to cut graphs into smaller pieces by
suitable boundary conditions in order to apply the ergodic theorem of
\cite{LenzMV}. These changes are finite rank perturbations in an appropriate sense
and hence change the spectral counting functions only in a very quantifiable
way, which is made precise in this section.

\medskip

Let ${H_1}$ and ${H_2}$ be selfadjoint
lower bounded operators with discrete spectra. Thus the
corresponding counting functions $n_{H_1}$ and $n_{H_2}$ are well-defined. Then,  the
spectral shift function is  the difference of two counting
functions, i.e.
\[
\xi_{{H_1},{H_2}}(\lambda) = n_{H_2}(\lambda) - n_{H_1}(\lambda).
\]

We will make use of the min-max principle in the form
$$
 n_H(\lambda) = \max \{ \dim X \mid X \leq \cD(H) , H|_X\leq \lambda \}
$$
where $X$ runs over all linear subspaces (denoted by $\leq$) of the domain of $H$ on which $H$ is bounded above
by $\lambda$ in the sense that $\langle f,Hf \rangle\leq \lambda \langle f,f\rangle$ for all $f\in X$.
The following lemma will be useful for dealing with several types of perturbations:
\begin{lem}\label{lem:SSFgeneral}
Let $H_1$, $H_2$ be selfadjoint lower bounded operators with discrete spectra on a separable Hilbert space.
Assume that $\cD_0\leq \cD(H_1) \cap \cD(H_2)$ such that $\cD_0$ has finite index $k$ in $\cD(H_2)$,
i.e.\ $\dim \cD(H_2)/\cD_0 =k$.
If $H_2|_{\cD_0} \geq H_1|_{\cD_0}$ then
$$  n_{H_2}(\lambda) \leq  n_{H_1}(\lambda) + k,\;\:\mbox{i.e.}\;\:
 \xi_{H_1,H_2}(\lambda) \leq k.  $$

\end{lem}
\begin{proof}
By the min-max principle for the counting functions we have:
\begin{align*}
n_{H_2}(\lambda) &= \max \{ \dim X \mid X\leq \cD(H_2),H_2|_X\leq \lambda \} \\
  &\leq \max \{ \dim X \mid X\leq \cD_0, H_2|_X\leq \lambda \} + k \\
  &\leq \max \{ \dim X \mid X\leq \cD_0, H_1|_X\leq \lambda \} + k \\
  &\leq \max \{ \dim X \mid X\leq \cD(H_1), H_1|_X\leq \lambda \} + k \\
  &= n_{H_1}(\lambda) + k
\end{align*}
and the proof is finished.
\end{proof}

This easily leads to the familiar boundedness of the spectral shift function for finite rank perturbations:
\begin{cor}
Let $H$ be selfadjoint with discrete spectrum and $K$ selfadjoint and bounded.
Let $Y_+$ (resp.\ $Y_-$) be the strictly positive (resp.\ negative) spectral subspace of $K$.
If $k_+:=\dim Y_+<\infty$ then $n(H,\lambda)-k_+\leq n({H+K},\lambda)$.
If $k_-:=\dim Y_-<\infty$ then $n({H+K},\lambda)\leq n(H,\lambda)+k_-$.
In particular, the spectral shift function for finite rank perturbations is bounded by the rank of the perturbation.
\end{cor}
\begin{proof}
Applying the previous lemma to $H_2=H+K$, $H_1=H$ (such that $\cD(H_2)=\cD(H_1)=\cD(H_1)$),
choosing $\cD_0=\cD(H)\cap Y_+^\perp$, yields one half of the claimed estimate.
The other half follows if we replace $K$ by $-K$ and $H$ by $H+K$.
\end{proof}
The perturbations in our applications will be of finite rank only in the following general sense:
 \begin{cor}\label{cor:SSFextensions}
Let $H_0$ be densely defined, closed, symmetric and lower bounded with deficiency index $k$.
Let $H_1$, $H_2$ be two selfadjoint extensions of $H_0$ with discrete spectrum. Then,
\begin{equation}
 |\xi_{H_1,H_2}| \leq k.
\end{equation}
\end{cor}
\begin{proof} This is a direct consequence of the previous lemma with $\cD_0 = D(H_0)$.
\end{proof}

\begin{rem}  Alternatively this corollary could be proven using
the invariance principle for the spectral shift function and the
second resolvent identity. In fact, this would not need the
assumption of discrete spectra.
\end{rem}

Let us now turn to the operators of interest for our applications.
These are the restrictions of $-\Delta$ to a finite subgraph of $G_d$, 
with selfadjoint  boundary conditions and 
possibly with an additional $L^\infty$-potential.

\begin{lem}\label{Randbedingungstoerung}
Let $G$ be a subgraph of $G_d$.  Let $\Delta$ and $\Delta'$ be two
  selfadjoint realizations of the Laplacian on $G$.  Set
  $\cD_0= \cD(-\Delta )\cap \cD(-\Delta')$ and $k:=\dim
  \cD(-\Delta)/\cD_0$.  Then $k\leq 2|E|$, and $\left|\xi_{-\Delta,
  -\Delta'} \right| \leq k $
\end{lem}
\begin{proof} Set $m:=|E|$. Obviously, 
$W_0^{2,2}(E):=\{f\in W^{2,2} (E) : (f(x),f'(x))=0 \ \forall x\in V(G)\}$ is
a subspace contained in all possible domains, and all of them are contained in $W^{2,2}(E)$. 
Because of $\dim W^{2,2}(E)/W_0^{2,2}(E)=4m$ and the rank $2m$ difference 
for the domains,
$\cD_0$ has index at most $2m$, and we can apply Lemma~\ref{cor:SSFextensions}.
\end{proof}

We are
particularly interested in Dirichlet resp.\ Neumann conditions on
the edges of a subgraph $G$.  The corresponding operators $-\Delta_D^G,-\Delta_N^G$
decompose completely as direct sums over the edges, and the counting
functions are given by $n_{-\Delta_D^G}(\lambda)=|E| n_D(\lambda)$ and  
$n_{-\Delta_N^G}(\lambda)=|E|n_N(\lambda)$.  Here,
\begin{align*}
n_D(\lambda) = \begin{cases} \lfloor \sqrt \lambda /\pi \rfloor &\text{if } \lambda\geq 0, \\ 0 &\text{otherwise,} \end{cases} \quad
n_N(\lambda) = n_D(\lambda) +\chi_{[0,\infty[}(\lambda)
\end{align*}
are the counting functions for the Dirichlet and Neumann Laplacians on $]0,1[$.
$\lfloor x\rfloor$ is defined for $x\in\RR$ by $\lfloor x\rfloor\in\ZZ$ and $x-1<\lfloor x\rfloor\leq x$.
 Together with the previous lemma, this implies  the following proposition.

\begin{prp}\label{akaRemark}   Let $G$ be a finite subgraph of $G_d$ and $\Delta$ any selfadjoint realization of the Laplacian on $G$. Then, 
 $\big| n_{-\Delta} (\lambda) - |E| n_D(\lambda) \big| \leq 2|E|$.
\end{prp}
In fact, looking at the quadratic forms one sees that the Dirichlet form dominates the forms of all other selfadjoint extensions;
therefore, $|E| n_D$ is a lower bound to the counting functions of all other selfadjoint extensions.

\medskip

Of course, everything said so far about selfadjoint extensions remains true when we add bounded perturbations.
The spectral shift is determined in the following lemma.
\begin{lem}  \label{Potentialstoerung}
Let $-\Delta$ be a selfadjoint realization of the Laplacian on a subgraph $G$ of $G_d$ and let  $W_1,W_2$ be bounded selfadjoint operators on $L^2(E)$.
Set $H_j=-\Delta +W_j,j=1,2$.
Then there is a constant $C$ depending only on the norms of $W_1,W_2$ such that
$$\left|\xi_{H_1,H_2} \right|(\lambda) \leq C |E|  $$
\end{lem}
\begin{proof}
Set $m=|E|$. According to Proposition \ref{akaRemark},
\[
\begin{array}{rcl}
 m\big(-2+ n_D(\lambda)\big) \leq &  n_{-\Delta_{}}(\lambda)  & \leq m\big( 2+   n_D(\lambda)\big) \quad\text{and therefore} \\
 m\big( -2+n_D(\lambda-\|W\|)\big) \leq & n_{-\Delta_{ }+W}(\lambda)  & \leq m\big(2+   n_D(\lambda+\|W\|)\big)
\end{array}
\]
for any bounded selfadjoint $W$ because $-\|W\|\leq W \leq \|W\|$ implies $n_H(\lambda-\|W\|)\leq n_{H+W}(\lambda)\leq n_H(\lambda+\|W\|)$ for any selfadjoint $H$.
Consequently, the explicit form of $n_D$  given above leads to
\[
\begin{array}{rcl}
 m\big( -3-\sqrt{\|W\|}/\pi \big) \leq & n_{-\Delta_{}+W}(\lambda) -\sqrt \lambda/\pi  & \leq m\big(2+ \sqrt{\|W\|}/\pi\big)
\end{array}
\]
and subtracting the estimates for $H_1$ and $H_2$ we finally get
\[
    \left|\xi_{H_1,H_2} \right|(\lambda) \leq m\big(5+ \sqrt{\|W_1\|}/\pi+\sqrt{\|W_2\|}/\pi\big)
\]
Using the fact that the Dirichlet counting function is a lower bound would allow us to replace $5$ by $3$.
\end{proof}

\begin{rem} The considerations of this section hold for general
finite metric graphs with all edges of equal length.
\end{rem}

\section{An abstract ergodic theorem}\label{s-ET}
In this section, we discuss the abstract ergodic theorem of \cite{LenzMV}
and show how it can be applied to give convergence of a suitably
defined spectral shift function of a family of random operators. 
See \cite{Lenz-02} for earlier one-dimensional results of the same type.
This is the crucial step in the proof of our main result.

\medskip

Let $\cala$ be a finite set and choose $\Gamma\in
\{\Gammat,\Gammaf\}$. Then, we can proceed as follows.  The set of all
finite subsets of $\ZZ^d$ is denoted by $\cals$.   A map $\varLambda \colon
\ZZ^d \longrightarrow \cala$ is called an \emph{$\cala$-colouring} of
$\ZZ^d$. A map $ P \colon Q(P) \longrightarrow \cala$ with $Q(P)\in
\cals$ is called an \emph{$\cala$-pattern}.  For $M\in \NN$ we denote by $C_M$ the cube at the origin with sidelength $M  - 1$, i.e.
$$C_M:=\{ x\in \ZZ^d:  0\leq x_j \leq M-1, j=1,\ldots, d\}.$$
The set of patterns $P$
with $Q(P) = C_M$ is denoted by $\calpM$. The set of all $C_M$, $M\in \NN$, is denoted by $\cC$, 
and a  pattern with  $Q(P) \in \cC$ is called a cube pattern.
For a pattern $P$ and $Q\in \cals$
with $Q \subset Q(P)$ we define the restriction $P\cap Q$ of $P$ to
$Q$ in the obvious way by $P\cap Q : Q \longrightarrow \cala$,
$x\mapsto P(x)$. For a pattern $P$ and $\gamma\in \Gamma$ we define
the shifted map $\gamma P$ by $ \gamma P : \gamma Q(P) \longrightarrow
\cala, \gamma (y) \mapsto P(y)$. On the set of all patterns we define
an equivalence relation by $P\sim P'$ if and only if there exists a
$\gamma\in\Gamma$ with $\gamma P = P'$. For a cube pattern  $P$ and an arbitrary pattern  $P'$ we define
the number of occurrences of the pattern $P$ in $P'$ by
\begin{equation*}
  \sharp^\Gamma_P P' :=
  \sharp\Bigl\{x\in Q(P') :
P'\cap \bigl(x +  Q(P)\bigr)\sim P\}.
\end{equation*}

\begin{dfn}
  \label{boundaryterm}
  A map $b \colon \cals \longrightarrow
  [0,\infty)$ is called a \emph{boundary term} if $b(Q) = b( t+ Q)$
  for all $t\in \ZZ^d$ and $Q\in \cals$, $\lim_{j\to \infty}
  |Q_j|^{-1} b (Q_j) =0$ for any van Hove sequence $(Q_j)$, and there
  exists $D>0$ with $b(Q) \leq D |Q|$ for all $Q\in \cals$.
\end{dfn}

\begin{dfn}
  \label{function}
  Let $(X,\|\cdot\|)$ be a Banach space and $F \colon
  \cals\longrightarrow X$ be given.    Let $\Gamma\in \{ \Gammat,\Gammaf\}$ be given.

  (a) The function $F$ is said to be \emph{almost-additive} if there
  exists a boundary term $b$ such that
  \begin{equation*}
    \left\| F ( \cup_{k=1}^m Q_k ) - \sum_{k=1}^m F(Q_k)\right\| \leq \sum_{k=1}^m b
    (Q_k)
  \end{equation*}
  for all $m\in\mathbb{N}$ and all pairwise disjoint sets $Q_k\in\cals$,
  $k=1,\ldots,m$.

  (b) Let $\varLambda :\ZZ^d \longrightarrow \cala$ be a colouring.
  The function $F$ is said to be \emph{$\Gamma$-$\varLambda$-invariant} if
  \begin{equation*}
    F( Q) = F (\gamma Q)
  \end{equation*}
  whenever $\gamma\in \Gamma $  and $Q\in\cals$ obey $\gamma (\varLambda
  \cap Q) = \varLambda \cap (\gamma Q)$.   In this case there exists a
  function $\widetilde{F}$ on the cubes $\cC$ with values in $X$ such that
  \begin{equation*}
    F( \gamma Q) = \widetilde{F} \Bigl(\gamma^{-1}\bigl(\varLambda\cap (\gamma  Q)\bigr)\Bigr)
  \end{equation*}
  for cubes $Q\in\cC $ and  $\gamma\in \Gamma$.

  (c) The function $F$ is said to be \emph{bounded} if there exists a
  finite constant $C>0$ such that
  \begin{equation*}
    \|F(Q)\| \leq C |Q|
  \end{equation*}
  for all $Q\in \cals$.
\end{dfn}

\begin{thm} \label{ergodictheorem}   Let $\cala$ be a finite set, $\varLambda \colon
  \ZZ^{d}\longrightarrow \cala$ an $\cala$-colouring and
  $(X,\|\cdot\|)$ a Banach-space.  Let $\Gamma \in \{\Gammaf,\Gammat\}$ be given. Let $(Q_j)_{j\in\mathbb{N}}$ be a van Hove
  sequence such that for every pattern  $P$ the frequency $\nu_P =
  \lim_{j\to \infty} |Q_j|^{-1} \sharp^\Gamma_P (\varLambda \cap Q_{j})$
  exists.       Let $F : \cals
  \longrightarrow X$ be a $\Gamma$-$\varLambda$-invariant, almost-additive
  bounded function. Then the limits
  \begin{equation*}
    \overline{F} := \lim_{j\to\infty} \frac{F(Q_j)}{|Q_{j}|}    = \lim_{M\to\infty}
    \sum_{P \in \calpM} \nu_P \frac{\widetilde{F}(P)}{|C_M|}
  \end{equation*}
  exist and are equal.
\end{thm}

\begin{rem}  (a) The theorem is proven in \cite{LenzMV} for $\Gamma = \Gammat$, However, the proof carries over verbatim to $\Gamma = \Gammaf$.

(b) We have explicit bounds on speed of convergence in terms of speed of convergence of the frequencies. For details we refer to \cite{LenzMV}.
\end{rem}

Now we turn to the case where the colouring is given by random potentials and boundary conditions.
  Then each $\omega\in \Omega$ induces a colouring
$\varLambda(\omega)$ with values in $\cA :=(\oplus_{j=1}^d \cB) \times \cL$ given
by
$$\varLambda(\omega) (x):=( (V([x,x+e_j])_j,L_x).$$

\begin{lem}  
Let $(Q_j)$ be an arbitrary van Hove sequence.  Then, for almost every
  $\omega\in \Omega$ the frequency $\nu_P = \lim_{j\to \infty}
  |Q_j|^{-1} \sharp^\Gamma_P (\varLambda (\omega) \cap Q_{j})$ exists and is independent of $\omega$ for every
 cube  pattern $P$.
\end{lem}
\begin{proof}
For a fixed pattern $P$ the frequency exists for almost every $\omega$
by a standard ergodic theorem. As there are only countably many $P$
the statement follows.
\end{proof}

In our setting we also have a family of random operators $(H_\omega)$
along with the Dirichlet Laplacian $\Delta_D$. As discussed above,
these operators can be restricted to the subgraphs $G_Q$ induced by
finite sets $Q$ of $\ZZ^d$. This yields the operator $H_\omega^Q$ and
$\Delta_D^Q$ with spectral counting functions $n_\omega^Q$ and $n_D^Q$
respectively. Now, $n_D^Q$ decomposes as a direct sum of
operators. Thus, denoting the eigenvalue counting function of the
negative Dirichlet Laplacian on $]0,1[$ by $n_D (\lambda)$ as above, we have
$n_D^Q = |E_Q| n_D = d |Q| n_D$.  The associated spectral shift
function is given as
$$
\xi_\omega^Q (\lambda) := n_\omega^Q (\lambda)- d \, |Q| \, n_D(\lambda) = d \, |Q| \,
\big(N_\omega^Q (\lambda)-n_D(\lambda)\big).
$$
The crucial point is that
$\xi_\omega$ falls into the framework of almost additive $F$
introduced above. This is shown next.

\begin{lem}
Let $(\cD,\|\cdot\|_\infty)$ be the Banach space of right
  continuous bounded functions on $\RR$. Then, for each $\omega\in
  \Omega$ the function $\xi_\omega : \cF \longrightarrow \cD$,
  $Q\mapsto \xi_\omega^Q$, is a bounded, $\varLambda(\omega)$ invariant
  almost additive function.
\end{lem}
\begin{proof} \textit{Almost additivity:}
Obviously, $Q\mapsto n_D^Q$ is almost additive with boundary term equal
  to zero, as Dirichlet boundary conditions decouple everything. As
  for $n_\omega$, we note that $ Q=\cup_{k=1}^m Q_k$ with pairwise
  disjoint $Q_k$, $k=1,\ldots,m$, induces a decomposition of $G_Q$
  into $G_{Q_k}$, $k=1,\ldots,m$. The restriction $H_\omega^Q$ of
  $H_\omega$ to $G_Q$ differs from the direct sum $\oplus_{k=1}^m
  H_{\omega^{Q_k}}$ only by boundary conditions on the boundary of the
  $G_{Q_k}$, $k=1,\ldots,m$. These boundary conditions yield 
  boundary terms in the corresponding eigenvalues counting functions
  by Lemma \ref{Randbedingungstoerung}. Thus, $n_\omega$ is almost
  additive as well.  Hence, $\xi_\omega$ is almost additive as it is a
  difference of almost additive functions.

\smallskip

\textit{$\varLambda(\omega)$-invariance}: This is clear from the definitions.

\smallskip

\textit{Boundedness:} By Lemma \ref{Randbedingungstoerung}, changing
all boundary conditions to Dirichlet conditions introduces an error
term of the order $2 |E|$. On the other hand, the operator with
potential and only Dirichlet conditions can easily be compared to the
Dirichlet operator without any potential by Lemma
\ref{Potentialstoerung}.

\end{proof}

\begin{rem}  The need to use a spectral shift function, i.e. the
difference between $n_\omega$ and $n_D$, in the above proof comes
exclusively from the last step in the proof. This use of a spectral
shift function to define the IDS bears some similarity to how the
integrated density of surface states is defined, see for instance
\cite{EnglischKSS-88,EnglischKSS-90,Chahrour-99a,Chahrour-99b,KostrykinS-01b}.
\end{rem}

The key result is now the following proposition.

\begin{prp} \label{conv} 
There is a bounded right continuous function $\Xi\colon \RR \to \RR$
such that for a given van Hove sequence $(Q_l)$ for  almost every $\omega \in \Omega$ the uniform convergence
\[
\lim_{l \to\infty}  \Big\|\frac{\xi_\omega^{Q_l}}{|E_{Q_l}|}-\Xi \Big\|_\infty =0\]
holds.
\end{prp}
\begin{proof}
Almost sure existence of the limit is a direct consequence of the
previous two lemmas and the first theorem of this section. In fact, this theorem gives an  explicit formula for the limit in terms of frequencies $\nu_P$ and
$\widetilde{F}$.  This formula shows  that the limit does not depend on $\omega$.
\end{proof}

\section{Proofs of the main results}\label{Proofs}
In this section we gather the material of the previous sections in order to
prove Theorem \ref{t-uniformIDS} and its two corollaries.

\bigskip

\begin{proof}[Proof of Theorem \ref{t-uniformIDS}:] We only consider the case $\Gamma = \Gammat$. The case $\Gamma = \Gammaf$ is similar and even simpler.

\smallskip

We first show independence of $N$ of the choice of $Q$:
Set $G_x=x+\bar S$. 
The invariance assumption on the random operators, the invariance of
the trace under unitary conjugation and the invariance of $\PP$ under
translations easily show that the expression
$$\int_\Omega \tr ( \chi_{G_{x}} \chi_{]\infty,\lambda]} (H_\omega))
d\PP(\omega)$$ does not depend on $x\in \ZZ^d$. Now, the claimed
independence of $Q$ follows easily.

\smallskip

We now show convergence of the $|E_{Q_l}|^{-1} n_\omega^{Q_l}$, $l\in \NN$. Proposition \ref{conv} yields almost sure convergence of
$$\frac{1}{|E_{Q_l}|} \xi_\omega^Q (\cdot) = \frac{1}{|E_{Q_l}|}
n_\omega^Q (\cdot )- n_D(\cdot)$$
to a limit $\Xi$ with respect to the supremum
norm. As the subtracted term $n_D(\cdot)$ does not depend on $l$ we
obtain  the desired convergence of $\frac{1}{|E_{Q_l}|} \xi_\omega^Q (\cdot) + n_D$.  The limit is $\widetilde{N}:= \Xi +n_D$.

\smallskip

Finally, we have to show that $N = \widetilde{N}$. To do so it
suffices to show that the measures associated to $N$ and
$\widetilde{N}$ respectively are equal. This in turn follows once we
show vague convergence. By a variant of the Stone/Weierstrass theorem it
suffices therefore to show that
$$ (*) \;\:\;\:\;\:\frac{1}{|E_{Q_l}|} \tr ( \chi_{G_{Q_l}} f(H_\omega) -
f(H_\omega^{Q_l}))\longrightarrow 0, \quad l\to \infty, $$
for all $f$ of the form $f(t)= (t- z)^{-1}$ for $z$ with non-vanishing
imaginary part. For $Q =Q_l$ we can split the graph $G_d$ into the
two components $G_Q$ and $G_d\setminus G_Q$.  Then, $H_\omega$ and
$H_\omega^{G_Q} \oplus H_\omega^{G_d\setminus G_Q}$ differ only by
boundary conditions on the set $V^\partial_Q$. Thus, by the second
resolvent identity

$$ D:=f(H_\omega) - f(H_\omega^{G_Q}\oplus H_\omega^{G_d\setminus
G_Q})$$ is an operator of rank at most $4 d |V^\partial_Q|$. Moreover,
$D$ is obviously bounded by $2 |Im (z)|^{-1}$. This gives
 \begin{align*} |\tr ( \chi_{G_{Q_l}} f(H_\omega) - f(H_\omega^{Q_l}))| = | \tr (
\chi_{G_{Q_l}} ( f(H_\omega)- f(H_\omega^{G_Q}\oplus
H_\omega^{G_d\setminus G_Q})))| \\
& \leq 8 d \frac{1}{|Im (z)|}
|V^\partial_{Q_l}|. \end{align*} 
As $(Q_l)$ is van Hove we obtain $(*)$.
\end{proof}

\begin{proof}[Proof of Corollary \ref{jump}:]
This is a variant of the proof of Theorem $2$ in \cite{KlassertLS-03}. 

$\Longrightarrow$: Let $\lambda$ be a point of discontinuity of $\mu$. Thus,
$\delta:=\mu(\{\lambda\})>0$. By the uniform convergence proven in Theorem
\ref{t-uniformIDS} this means that the multiplicity of the eigenvalue
$\lambda$ of $H_\omega^{Q_l}$ is at least $\delta/2 |E_{Q_l}|$ for $l$
large.  On the other hand,    within the eigenspace of $H_\omega^{Q_l}$ to the eigenvalue $\lambda$ the space of 
 $f$ with $(f(v),f'(v))=0$ for all
$v\in V^\partial_{Q_l}$ has codimension at most $2 d
|V^\partial_{Q_l}|$. As $(Q_l)$ is van Hove, we have$$ 2 d
|V^\partial_{Q_l}|< \delta/2 |E_{Q_l}|$$ for $l$ large enough and we
obtain a compactly supported eigenfunction of
$H_\omega^{Q_l}$ to $\lambda$ all of whose boundary values are $0$. Thus, this
function can be extended by $0$ to give a compactly supported
eigenfunction.

\smallskip

$\Longleftarrow$: Choose the van Hove sequence to consist of cubes centered around the origin with increasing sidelengths. For $l\in \NN$, define the function $\chi_l$ on
$\Omega$ by $\chi_l (\omega)=1$ if there exists a compactly  supported
eigenfunction of $H_\omega$ to $\lambda$ which is supported in $Q_l$ and
vanishes identically on all edges $e\in E^\partial_{Q_l}$ and $\chi_l
(\omega) =0$ otherwise. Then
$$0< \delta_l:= \int_\Omega \chi_l (\omega)d\PP(\omega)$$
for all $l$ large enough. Fix such an $l$ and set $\chi:=\chi_l$, $\delta:=\delta_l$ and $R$ to be the diameter of $Q_l$.  Then, the ergodic theorem yields
$$\delta=\lim_{k\to \infty} \frac{1}{|Q_k|} \sum_{x\in Q_k} \chi(T^x \omega)$$
almost surely.  Choose $k_0$ with $\frac{1}{|Q_k|} \sum_{x\in Q_k} \chi(T^x \omega)\geq \delta/2$ for $k\geq k_0$.
For each $k\geq k_0$ we consider a maximal system $M_K$ of points in $Q_k$ such that $\chi(T^x\omega)=1$ and points in $M_k$ have distance at least $2R$. A
 direct combinatorial  argument then shows that
$$|M_k|\geq \frac{\delta |Q_k|}{2 |B_{2R}|}$$
By construction different points in $M_k$ yield orthogonal compactly supported eigenfunctions of $H_\omega$.
This yields $$\mu_\omega^{Q_k} (\{\lambda\}) \geq \frac{\delta |Q_k|}{2 |B_{2R}|}$$
for $k>k_0$ and the statement follows from the uniform convergence.
\end{proof}

\begin{proof}[Proof of Corollary \ref{spectrum}:]
Due to ergodicity and the Pastur-Shubin trace formula the corollary follows by standard arguments.

\end{proof}

\section{Operators with magnetic fields}\label{magnetic}
Our setup is general enough to include magnetic fields as well. Since
this is not quite obvious from the definitions and since there is a
special interest in these operators we devote a section to the question
how magnetic fields enter through boundary conditions.

To this end, let $G$ be a metric graph and $L$ a choice of boundary conditions
as in Sections \ref{s-MR} and \ref{s-FRP}. The most general symmetric
first order perturbation of $-\frac{\dd^2}{\dd t^2} $ on an edge $e\in
E$ is, up to zeroth order terms, given by 
\[ H(a)_e := -\left(\frac{\dd}{\dd t}-\imath a_e\right)^2 \]
for arbitrary real valued $a_e\in C^1(\bar e)$, where $\bar e$ is the
closure of the edge $e$, i.e.\ identified with the closed interval
$[0,1]$. If $G$ is embedded, e.g.\ in $\RR^d$ such as
in our case, it is natural to fix a so called vector potential $A\in
\left((C^1(\RR^d)\right)^d$ and to define $a_e$ to be the component
$A_j$ of $A$ if $e$ is embedded with euclidean direction $e_j$. This
corresponds to a magnetic field $2$-form $B=\sum_{ij}\frac{\partial
A_j}{\partial x_i}\,\dd x_i\wedge \dd x_j$ on $\RR^d$. We will work in
the general setting, though.

The selfadjoint realization of $H(a)$ corresponding to $L$ is then given
by the domain
\[ \cD(H_L(a)) = \{f\in W^{2,2}(E): \forall x\in V: (f(x),f'(x)-\imath
(af)(x)) \in L_x \} \]
as the usual partial integration argument shows;
i.e.\ one has to specify mixed Dirichlet and (magnetic) Neumann boundary
conditions as expected.

Now define the phase along the edge $e$ by $\varphi_e(t)=\int_0^t
a_e(s)\,\dd s$. Since each $a_e$ is real valued, multiplication by
$e^{\imath \varphi_e}$ on each edge $e$ defines a unitary operator $U$ on
$L^2(E)$ which maps $W^{2,2}(E)$ into itself. A simple calculation shows that
\[ \left( U^* H(a) U\right)_e = -\frac{\dd^2}{\dd t^2} \]
on $W^{2,2}(E)$. Therefore, $H_L(a)_e$ is unitarily equivalent to
$H_e=-\frac{\dd^2}{\dd t^2}$ with domain
\begin{align*} \cD(H_{\tilde L}) &= \{f\in W^{2,2}(E): Uf\in \cD(H_L(a)) \}
\\
 &= \{f\in W^{2,2}(E): \forall x\in V: (( Uf)(x),(Uf)'(x)-\imath
(aUf)(x)) \in L_x \} 
\end{align*}
for some graph local boundary condition $\tilde L$ since multiplication
by $U$ preserves graph locality. In order to determine $\tilde L$ we
calculate, using the notation from Section~\ref{s-MR}:
\begin{align*}
 (Uf)_e(s(e)) &= f_e(s(e)) \\
 (Uf)_e(r(e)) &= e^{\imath \varphi(r(e))}f_e(r(e)) \\
 (Uf)'_e(s(e))-\imath(aUf)(s(e)) &= f'_e(s(e)) \\
 (Uf)'_e(r(e))-\imath(aUf)(r(e)) &= e^{\imath \varphi(r(e))} f'_e(r(e))
\end{align*}
Hence we define a unitary diagonal matrix $u$ on $\CC^{4d}$ by setting
a diagonal entry to $e^{\imath \varphi(s(e))}=1$ if it corresponds to a boundary value at some
$s(e)$ and to $e^{\imath \varphi(r(e))}$ if it corresponds to a boundary
value at some $r(e)$. Then $\tilde L = u^* L$.

Note that Dirichlet and (magnetic) Neumann conditions for the
magnetic operator $H(a)$ transform under $u$ into Dirichlet and
(non-magnetic) Neumann conditions for $H$. This is to be expected because
these boundary conditions decouple, and any magnetic field can be gauged
away in strictly onedimensional systems. In contrast, (magnetic)
Kirchhoff conditions for $H(a)$ transform into completely different
boundary conditions: At a vertex $x$, only the values at outgoing edges
(those with $s(e)=x$) coincide; they, in turn, coincide with the values
$e^{\imath \varphi(r(e))}f_e(r(e))$ at the incoming edges so that $f$
need not be continuous at the vertex. Similarly, the derivative part of
the Kirchhoff condition is modified to a weighted condition.

Finally, let us note that the results of this section imply that the
results of section~\ref{s-MR} for Schr\"odinger operators with random
potentials and random (or fixed) boundary conditions lead to the same
results for magnetic Schr\"odinger operators with random potentials and
random (or fixed) magnetic fields, specified through the phases
$\varphi_e(r(e))$ at the endpoints.

\section{Percolation on $G_d$}\label{Site}
In this section we discuss two models which can be seen as versions of
site and edge percolation on $G_d$, respectively. A third model corresponds to 
percolation on the junctions between vertices and edges. These models are
based solely on random boundary conditions. The potential of the operators is identically equal to zero.
 Unlike in the percolation models on combinatorial graphs ``deleted'' edges are not
removed completely from the graph but only cut off by Dirichlet boundary
conditions. The reason is that removing edges would mean to remove
infinite dimensional subspaces from our Hilbert space. This would result 
in a spectral distribution function which is not comparable to the one of
the Laplacian with the concerned edge included.  Hamilton operators on percolation subgraphs of combinatorial graphs
have been considered in the literature in theoretical physics \cite{deGennesLM-59a,KirkpatrickE-72,ChayesCFST-86},
computational physics  \cite{KantelhardtB-02} (and references therein), and mathematical physics
\cite{Veselic-05a,Veselic-05b,KirschM-06,MuellerS,Veselic-06}. 

Before giving details we would like to emphasize the following: The
examples below include cases in which the graphs  contain many 
finite components giving rise to compactly supported eigenfunctions.
In particular,
the integrated density of states has many discontinuities. 
In fact, in the subcritial phase the IDS is a step function.
However, despite all these jumps our result
on uniform convergence does hold!

\bigskip

We first discuss  a site percolation model. The percolation process is defined by the following procedure:
toss a (possibly biased) coin at each vertex and
-- according to the outcome -- put either a Dirichlet or a Kirchhoff
boundary condition on this vertex. Do this at every vertex independently of all the others.
 To be more precise, let $p\in (0,1)$ and $q=1-p$ be given. 
Let $\cA :=\{L^D, L^K\}$ and the probability measure $\nu:= p \delta_{L^K}+ (1-p)\delta_{L^D}$ on $\cA$ be given. 
Define $\Omega$ as the cartesian product space $\times_{x \in V_d} \cA$ with product measure
$\PP:=\otimes_{x\in V_d} \nu$. Let $L$ be the stochastic process 
with coordinate maps 
$L(\omega)(x) := \omega (x)$. These data yield a family of random operators 
$-\Delta_\omega :=H_\omega$  acting like the free Laplacian with domain given by
$$D( \Delta_\omega)=\{f\in W^{2,2} (E) : (f(x),f'(x))\in
\omega(x)  \ \forall x \in V_d\}.$$
Intuitively, placing a Dirichlet boundary condition at a vertex means ``removing'' it from the
metric graph. The $2d$ formerly adjacent edges have now ``loose ends''.
A fundamental result of percolation theory tells us that 
for sufficiently small values of $p$ the percolation graph consists entirely of finite components almost 
surely.
For these values of $p$ our Laplace operators decouple
completely into sums of operators of the form $-\Delta^G$ for finite
connected subgraphs $G$ of $G_d$. Here, $\Delta^G$ acts like the free
Laplacian and has  Dirichlet boundary conditions on its deleted vertices (boundary vertices) and 
Kirchhoff boundary conditions in its vertices which have not been deleted by the percolation process 
(interior vertices). 
We introduce an equivalence relation on the set of connected subgraphs of 
$G_d$ with a finite number of edges by setting $G^1 \sim G^2$ iff there 
exists 
a $\gamma \in \Gamma$ such that $ \gamma G^1 = G^2$. For such an 
equivalence class $\cG$ we define 
$n^\cG$  as the eigenvalue counting function of $-\Delta^G$ for some $G \in 
\cG$, 
and set $N^\cG= \frac{n^\cG}{|E_G|}$.
Defining the density  $\nu_\cG$
of an equivalence class of  finite subgraphs of $G_d$ within the 
configuration $\omega$ in the obvious way, 
we obtain
as integrated density of states for the family $H_\omega$
$$ N = \sum_{\cG} \nu_\cG N^\cG,$$ where the sum runs over all equivalence 
classes $\cG$ of finite
connected subgraphs of $G_d$. 
Thus, the integrated density of states is a pure
point measure in this case with many jumps. More interestingly, all
these jumps remain present (even if their height is diminished) when we
start increasing $p$. This yields models in which the operators are
not given as a direct sum of finite graph operators but still have
lots of jumps in their integrated density of states.
Related phenomena for combinatorial Laplacians
have been studied e.g.~in \cite{ChayesCFST-86,Veselic-05b,Veselic-06}.

\bigskip

We now discuss an edge percolation model.  The basic idea is to
decide for each edge independently whether Dirichlet boundary conditions
are put on both ends or not. All other boundary conditions are
Kirchhoff type. The problem when defining this edge percolation model is that our stochastic processes
are indexed by vertices rather than edges. We thus have to relate
edges to vertices. This is done by going to each vertex and then tossing
a (biased) coin for each $j=1,\ldots, d$ to decide how to deal with
the edge $[x,x+e_j]$.

Here are the details. Let $p_0\in (0,1)$ and
$p_1= 1 - p_0$ be given.  Let $\cA$ consist of all maps $S$ from
$\{1,\dots,d\}$ to $\{0,1\}$.  Put a probability measure $\nu$ on
$\cA$ by associating the value $\prod_{j=1}^d p_{S(j)}$ to the element
$S$. Now, $\Omega$ is the cartesian product space $\times_{x \in V_d} \cA$  with product
measure $\PP:=\otimes_{x\in V_d} \nu$. To each $\omega\in \Omega$ we
associate the operator $-\Delta_\omega =H_\omega$ which acts like the
free Laplacian and has boundary conditions as follows: The edge $e
=[x,x+e_j]$ has Dirichlet boundary conditions on both ends if the random variable associated to the vertex $x$
has the $j$-th component equal to $1$.  Otherwise the boundary condition is chosen to be
Kirchhoff. Here, again the operator decouples completely into 
operators on finite clusters for small	
enough values of $p_0$.

\bigskip

Similary one can consider a percolation process indexed by pairs $(x,e)$
of adjacent vertices and edges. As in the last model consider a colouring
$\cA$ consisting of all maps $S$ from $\{1,-1,2,-2\dots,d, -d\}$ to $\{0,1\}$.
The probability space and measure are defined similarly as before.
Each $\omega$ gives rise to a Laplace operator with the following boundary conditions:
if the $-j$-th component of the random variable associated to the vertex  $x$ has the value one,
then the edge $[x-e_j,x]$ is decoupled from $x$ by a Dirichlet boundary condition.
If the $j$-th component of the same random variable has value one then the edge $[x,x+e_j]$ is decoupled from $x$ by a Dirichlet boundary condition. Conversely, those components of the random variable which are zero
correspond to Kirchhoff boundary conditions.

\vspace*{2cm}

{\small\textbf{Acknowledgment}
It is a pleasure to thank Mario Helm for interesting discussions.
The authors were financially supported by the DFG,  two of them (M.~G.~and I.V.) 
under grant   Ve 253/2-2 within the Emmy-Noether-Programme.}
 
 \bibliographystyle{alpha}
 \bibliography{quantum-g}

\end{document}